\documentclass[a4paper,12pt,oneside]{article}
\usepackage{graphicx}
\usepackage[all]{xy}
\usepackage[centertags]{amsmath}
\usepackage{latexsym}

\usepackage{amsfonts}

\usepackage{amssymb,amsthm}
\usepackage[english]{babel}

\usepackage{newlfont}

\usepackage{indentfirst}
\usepackage[english]{babel}
\usepackage[latin1]{inputenc}
\usepackage{eufrak}
\newtheorem{remark}{Remark}[section]

\theoremstyle{plain}
\newtheorem{lemma}{Lemma}[section]
\newtheorem{proposition}{Proposition}[section]
\newtheorem{theorem}{Theorem}[section]
\newtheorem{definition}{Definition}[section]
\newtheorem{corollary}{Corollary}[section]

\newtheorem{example}{Example}[section]

\newcommand{\R}{\mathbb R}
\newcommand{\C}{\mathbb C}

\newcommand{\G}{\mathcal{G}}%%%

\newcommand{\beqn}{\begin{eqnarray}}
\newcommand{\eeqn}{\end{eqnarray}}
\newcommand{\beq}{\begin{eqnarray}}
\newcommand{\eeq}{\end{eqnarray}}
\newcommand{\bpro}{\begin{proposition}}
\newcommand{\epro}{\end{proposition}}
\newcommand{\blem}{\begin{lemma}}
\newcommand{\elem}{\end{lemma}}
\newcommand{\bdfn}{\begin{definition}}
\newcommand{\edfn}{\end{definition}}
\newcommand{\bcor}{\begin{corollary}}
\newcommand{\ecor}{\end{corollary}}
\newcommand{\bthm}{\begin{theorem}}
\newcommand{\ethm}{\end{theorem}}
\newcommand{\bex}{\begin{example}}
\newcommand{\eex}{\end{example}}
\newcommand{\brmq}{\begin{remark}}
\newcommand{\ermq}{\end{remark}}
\newcommand{\benum}{\begin{enumerate}}
\newcommand{\eenum}{\end{enumerate}}
\newcommand{\bitem}{\begin{itemize}}
\newcommand{\eitem}{\end{itemize}}

\theoremstyle{plain}

\usepackage[latin1]{inputenc}

\title{On properties of principal elements of Frobenius Lie algebras}

\author{ Andre Diatta  \footnote{ \footnotesize    CNRS Institut Fresnel UMR  7249,  Aix-Marseille Universit\'e, Marseille, France
\newline
E-mail: andre.diatta@fresnel.fr
}
 and Bakary Manga\footnote{ \footnotesize
(2) D\'epartement de Math\'ematiques et Informatique. Universit\'e Cheikh Anta Diop de Dakar, S\'en\'egal. 
\newline
E-mail: bakary.manga@ucad.edu.sn
}
 }

\begin{document}

\maketitle

%\author{Andr\'e Diatta \\
%University of Liverpool. Department of Mathematical Sciences,
%\\
 %Peach Street, Liverpool L69 7ZL, UK}

\begin{abstract}
 We investigate  the properties of principal elements of Frobenius Lie algebras.
 We prove that any Lie algebra with a left symmetric algebra structure can be embedded 
as a subalgebra of some $\EuFrak{sl}(m,\mathbb K),$ for $\mathbb K=\mathbb R$ or $\mathbb C$. Hence, the work of Belavin and Drinfeld on solutions of the Classical Yang-Baxter Equation on simple Lie algebras, applied to the particular case of $\EuFrak{sl}(m,\mathbb K)$ alone, paves the way to the complete classification of Frobenius and more generally quasi-Frobenius Lie algebras.
We prove that, if a Frobenius Lie algebra has the property that every derivation is an inner derivation, then every principal element is semisimple.
 As an important case, we prove that in the Lie algebra of the group of affine motions of the Euclidean space of finite dimension, every derivation is inner.
We also bring examples of Frobenius Lie algebras that are subalgebras of $\EuFrak{sl}(m,\mathbb K),$  
but yet have  nonsemisimple principal elements as well as some with semisimple principal elements  
having nonrational eigenvalues, where $\mathbb K=\R$ or $\mathbb C.$
\end{abstract}

%%====================

\section{Introduction }

Throughout this work, all Lie algebras are considered over $\mathbb K=\mathbb C$ or $\mathbb R.$ 
A finite dimensional Lie algebra $\G$ with Lie bracket $[.,.]$ over $\mathbb K,$ is 
called quasi-Frobenius or symplectic, if it possesses a nondegenerate skew-symmetric
 closed 2-form $\omega.$ If $\omega$ is the (Chevalley-Eilenberg) differential 
 of some linear functional, i.e. if  there exists a linear 1-form $\alpha$ on $\G,$
  such that $\omega (x,y)= \partial \alpha (x,y) := - \alpha ([x,y]),$  for all
 $x,y$ on $\G,$ then $\G$ is called a Frobenius (or exact symplectic) Lie algebra 
and $\alpha$ a Frobenius functional. Quasi-Frobenius Lie algebras are also endowed with a natural left symmetric algebra (LSA) structure (J.L. Koszul, A. Joseph, A. Albert, Y. Matsushima, ... )  See Section \ref{preliminaries} for more details.
Frobenius Lie algebras have been extensively studied by several authors, see e.g. \cite{bordemannmedinaouadfel}, \cite{bon-yao-chu}, \cite{dardie-medina-symplectic}, \cite{dardie-medina-kahler}, \cite{dergachev-kirillov}, \cite{Do-Na}, \cite{Elashvili}, \cite{Gerstenhaber97}, \cite{Gerstenhaber2009}, \cite{Gerstenhaber2008-arxiv}, \cite{Gindikin-Sapiro-Vinberg}, \cite{Li-Me},  \cite{Lichne91},  \cite{ooms74},  \cite{ooms80}, \cite{Rais}, \cite{Shapiro69}, etc.
  They also
 appear in many research areas such as Poisson Geometry and Hamiltonian systems (see e.g. \cite{Dr83}, \cite{Diatta2000}),  invariant affine geometry (see. e.g. \cite{di-me-cybe}), symplectic and K\"ahler geometry, homogeneous domains (\cite{bordemannmedinaouadfel}, \cite{dardie-medina-symplectic}, \cite{dardie-medina-kahler}, \cite{Do-Na}, \cite{Li-Me}, \cite{Gindikin-Sapiro-Vinberg}, \cite{Shapiro69}), the Yang-Baxter Equation (see e.g.  \cite{belavindrinfeld82, belavindrinfeld83, belavindrinfeld98}, \cite{stolin}),  contact geometry   (\cite{Diatta2000}, \cite{diatta-contact}), deformation quantization, etc.

In \cite{Gerstenhaber2009}, M. Gerstenhaber and A. Giaquinto, proved that if a Frobenius Lie algebra $\mathcal G$ is a saturated subalgebra (i.e. $\mathcal G$  is not an ideal of any larger subalgebra)  of a simple Lie algebra, then any principal element of a Frobenius functional is semisimple, at least for $\mathbb K$ algebraically closed.
In this paper, we drop the saturation requirement and explore the properties of principal element of Frobenius Lie algebras, in general.
\\
$\bullet$  We explore the geometric properties of principal elements of general Frobenius Lie algebras, from the viewpoint of invariant affine geometry, namely, by exploiting the induced structure of left symmetric algebra (Section \ref{chap:properties_of_principal_and_rightnil}).
\\
$\bullet$ In Theorem \ref{affineembedding}, we prove that, more generally, any Lie algebra with a left symmetric algebra structure (hence including quasi-Frobenius and in particular Frobenius Lie algebras) can be embedded as a subalgebra of some $\EuFrak{sl}(m,\mathbb K).$ 
This implies that some of the results within \cite{Gerstenhaber2009}, \cite{Gerstenhaber2008-arxiv}, ...,  are valid for a  much wider class of Frobenius Lie algebras.
In particular, for any Frobenius Lie algebra, one can still perform the study on associated graphs as in, e.g. \cite{Gerstenhaber2008-arxiv}, amongst other papers. Such an embedding of a Frobenius Lie algebra in some $\EuFrak{sl}(m,\mathbb K),$   is not unique. It would be interesting to investigate the obstructions to embedding a given Frobenius Lie algebra as a subalgebra of  $\EuFrak{sl}(m,\mathbb K),$ of some special given type, e.g. seaweed type, etc.
Another consequence is that, the classification of  quasi-Frobenius  (and Frobenius) Lie algebras of dimension $2n,$ could solely rely on that of the solutions of the Classical Yang-Baxter Equation on  $\EuFrak{sl}(2n+1,\mathbb K).$ Amongst other authors,
Belavin and Drinfeld  \cite{belavindrinfeld82, belavindrinfeld83, belavindrinfeld98}  have dedicated some of their work to the classification of the latter on simple Lie algebras. However, the classification of Frobenius Lie algebras might need to take into account the so-called double extension process  as in \cite{dardie-medina-symplectic}, \cite{dardie-medina-kahler}.
\\
$\bullet$
We prove that, if a Frobenius Lie algebra has the property that every derivation is an inner derivation, then every principal element is semisimple (Theorem \ref{thm:semisimplecohomology}), at least for $\mathbb K=\mathbb C.$  As an important case, we prove that in the Lie algebra of the group of affine motions of the Euclidean space of finite dimension, every derivation is inner (Theorem \ref{thm:affineliealgebra}).
\\
$\bullet$  In Section \ref{chap:examples}, we also supply a  wide family of interesting examples of Frobenius Lie algebras, that hence  can be considered as subalgebras of some $\EuFrak{sl}(m,\mathbb K),$ the adjoint of whose principal elements could be nearly anything. In particular,  some of them have  nonsemisimple principal elements, some other have semisimple principal elements  the adjoint of which have nonrational eigenvalues, ..., for $\mathbb K=\R$ or $\mathbb C.$

\medskip

The image $Im(r)$ of a solution $r$ of the  Classical Yang-Baxter Equation on a Lie algebra $\G,$ is a quasi-Frobenius Lie algebra which is a subalgebra of $\G$ (see e.g. \cite{Li-Me}).  If $\mathcal G$ is  compact, then $Im (r)$ is an abelian subalgebra (see e.g. \cite{Li-Me}). Hence, a Frobenius Lie algebra cannot be a subalgebra of a compact Lie algebra.
 On the other hand, if a quasi-Frobenius Lie algebra $\G$ is unimodular (see e.g. \cite{bon-yao-chu}), or equivalently, if the connection (\ref{eq:LSA-symplectic}) is complete on any connected Lie group with Lie algebra $\G$ (\cite{Li-Me}), then $\G$ is solvable. If in addition, it is K\"ahler, then it is Abelian (\cite{Li-Me}). Of course, the results within this paper clearly indicate that, although every Frobenius Lie algebra is a subalgebra of some $\EuFrak{sl}(n,\mathbb K),$ it needs not be solvable, although it is obviously the image of a solution of the Classical Yang-Baxter Equation.

\section{Some preliminaries and notations}\label{preliminaries}
 Throughout this paper, we will always denote by $\langle ,\rangle$ the duality pairing between a vector space $E$ and its dual $E^*$, unless otherwise stated. If $(e_1,\ldots, e_n)$ is a basis of $E,$ then the corresponding dual basis will be often denoted by $(e^*_1,\ldots, e_n^*).$
 Let $G$ be a Lie group with unit $\epsilon$  and $\G$ its Lie algebra identified with the tangent space $T_{\epsilon}G$ to $G$ at $\epsilon.$ If $x$ is an element of $\G,$ we denote by $x^+,$ the unique left invariant vector field on $G,$ with value $x=x^+_{\epsilon}$ at $\epsilon.$ That is, for every $\sigma\in G,$  if we let $L_\sigma$ stand for the left multiplication $\tau\mapsto L_{\sigma}(\tau):=\sigma\tau$  on $G$ by $\sigma,$  and $T_{\epsilon}L_{\sigma}$ its derivative at $\epsilon,$ then $x^+_{\sigma}:=T_{\epsilon}L_{\sigma}x.$  For all the representations involved in this work, we will equally use the same symbol $\partial$ for the associated 
(Chevalley-Eilenberg) coboundary operator. 

\begin{definition}
 A pair $(G,\omega^+)$ of a Lie group $G$ and a left invariant nondegenerate closed differential 2-form
 $\omega^+$ on $G,$ is termed a quasi-Frobenius (or symplectic) Lie group. If $\omega:=\omega^+_{\epsilon},$ then $(\mathcal G,\omega)$ is the corresponding quasi-Frobenius (or symplectic) Lie algebra.
   If in addition, $\omega^+$ is the differential  $d \alpha^+=\omega^+$ of a left invariant differential 1-form $\alpha^+$ on $G$, or equivalently, if $\omega(x,y)=\partial \alpha(x,y):=-\alpha ([x,y])$ for all $x,y\in\G,$ with $\alpha\in\G^*,$ then   $(G, \omega^+)$ (resp. $(\mathcal G,\omega)$ ) is called a Frobenius (or exact symplectic) Lie group (resp. Lie algebra.)
 \end{definition}
   From now on, we may as well use the duality pairing notation $\langle\alpha, [x,y]\rangle$ 
instead of $\alpha ([x,y]).$
  Now consider the linear map
 \beqn\label{eq:q}
 q:\mathcal G \to \mathcal G^*, ~~~ q(x):=i_x\omega=\omega(x,.).
 \eeqn
  Then  $\omega$ being a closed 2-form is equivalent to $q$ being a cocycle for the co-adjoint action, 
 that is $0=\partial q (x,y):=ad_x^*q(y)-ad^*_yq(x)-q([x,y]),$ for all $x,y\in\G.$ Moreover, $\omega$ is the 2-coboundary $\omega=\partial \alpha$ of $\alpha\in\G^*$ for the adjoint action, if and only if $q$ is the 1-coboundary of $-\alpha$ for the coadjoint action, namely,
$q(x)=ad_x^*\alpha,$ for every $x\in \mathcal G.$ Of course, the nondegeneracy  of $\omega$ is equivalent to the invertibility of $q.$

Recall that a manifold $M$ is referred to as an affine manifold,
if it is endowed with an affine structure.
That is, if $M$ has a maximal atlas of charts (affine atlas)
 to an affine space (say $\R^n$) with transition functions
extending to affine transformations. Recall also that, affine transformations
of $\R^n$ form the group
of transformations that preserve the usual
 flat and torsion free
connection on $\R^n.$
Hence, every affine chart inherits a
 flat and torsion free connection, the pullback of the usual one of $\R^n.$ As the change of coordinates are affine, those connections glue together into a globally defined torsion free and (locally) flat connection, say $\nabla,$ on $M.$
Conversely, any torsion free and (locally) flat connection on $M$ gives rise to an affine structure on $M.$
A vector field $X$ on an affine manifold $(M,\nabla)$ is said to be radiant
if it satisfies $\nabla_YX=Y$ for every vector field $Y$
 on $M.$ A map between two affine manifolds whose expression
in affine atlas, extends to an affine map,
 is called a morphism of affine manifolds.
An affine structure on a Lie group $G$ is said to be left invariant if left translations
 are morphisms of the underlying affine manifold. The restriction of the corresponding connection to left invariant vector fields gives rise to
a product $(x,y)\mapsto xy$ on the Lie algebra $\G$ of $G$ by
\beqn \label{eq:LSA-of-connection}
xy:=(\nabla_{x^+}y^+)_{\epsilon}.
\eeqn
If we let $A(x,y,z)$ stand for the associator $A(x,y,z):= (xy)z-x(yz)$ of $x,y,z\in\G$,
then the (local) flatness of $\nabla$ is then equivalent to $A$ being left symmetric,
i.e.
\beqn \label{eq:LSA-associator}
A(x,y,z)=A(y,x,z), \text{~i.e.,~} (xy)z-x(yz)=(yx)z-y(xz), \text{$~\forall$} x,y,z\in \G.
\eeqn
The torsion free condition is equivalent to the fact that the Lie bracket $[x,y]$ of any two vectors $x,y$ is their commutator
\beqn \label{eq:LSA-torsion}
xy-yx=[x,y].
\eeqn

\begin{definition}
A Lie algebra $\G$ with a product $(x,y)\mapsto xy,$ satisfying the conditions (\ref{eq:LSA-associator}) and (\ref{eq:LSA-torsion})
is called a left symmetric algebra (LSA, for short). Such a product is called an LSA structure on $\G.$
Sometimes, we will abusively refer to such  a product as an LSA.
\end{definition}

Let us define  a linear map  $\mathcal L$ from $\G$ to the space $\EuFrak{g l}(\G)$ of all linear
maps from $\G$ to itself, as follows:
$\mathcal L_x(y):=xy$.
 The condition (\ref{eq:LSA-associator}) above is equivalent to $\mathcal L$ being a homomorphism of Lie algebras and hence
a representation of $\G$ on itself, where the Lie algebra structure of $\EuFrak{g l}(\G)$ is given by the matrix commutators. The corresponding $\G$-module will be
 denoted by $\G_{\mathcal L}.$
From that point of view, the condition (\ref{eq:LSA-torsion}) is equivalent to the identity transformation $\mathbb I(x):=x,$
from $\G$ to $\G_{\mathcal L},$ being  a 1-cocycle for the representation $\mathcal L.$ That is,
$0=\partial \mathbb I(x,y):=\mathcal L_x\mathbb I(y)-\mathcal L_y\mathbb I(x)-\mathbb I([x,y]).$
Conversely, associated to an LSA structure on a Lie algebra $\G,$ is a left invariant affine structure on any Lie group with Lie algebra $\G.$ The corresponding  left invariant flat and torsion free connection  is obtained by left translating about the group, the products $xy$ of all $x,y\in\G$:
\beqn\label{eq:connection-of-LSA}
\nabla_{x^+}y^+:=(xy)^+.
\eeqn
Hence, (\ref{eq:LSA-of-connection}) and (\ref{eq:connection-of-LSA}) provide a bijective correspondence between left invariant affine structures in connected Lie groups and  left symmetric algebra (LSA) structures in their Lie algebras.

\begin{definition} Let $\G$ be a Lie algebra endowed with  an LSA $(x,y)\mapsto  xy$.
(1) We call $y_0$ a right-unit of the LSA, if the right multiplication $\mathcal R_{y_0}: x\mapsto xy_0$ by $y_0,$ coincides with the identity map $\mathbb I_{\mathcal \G}$  of $\G.$ That is $x=xy_0,$ for all $x\in\G.$
(2) An element $\tilde x_0$ is termed a right-nil of the LSA, if the right multiplication $\mathcal R_{\tilde x_0}$ by $\tilde x_0,$ coincides with the zero map  of $\G.$ That is $0=y\tilde x_0=\mathcal R_{\tilde x_0}(y),$ for all $y\in\G.$
\end{definition}

A left invariant vector field $x^+_0$ on a Lie group $G$ is radiant if and only if, for the corresponding LSA on $\G,$  the vector $x_0:=(x^+_0)_{\epsilon}$ is a right-unit.
  The existence of a right-unit $x_0$  can be understood as the fact that, for the Lie algebra representation $\mathcal L,$ the 1-cocycle $\mathbb I$ is actually equal to the coboundary $\partial x_0$ of $x_0\in\G_{\mathcal L},$ namely,  $y=\mathbb I(y)=\partial x_0 (y):=\mathcal L(y)x_0:=yx_0,$ for all $y\in\G.$

To any quasi-Frobenius Lie group $(G,\omega^+)$ is associated a zero curvature and torsion free linear  connection $\nabla$ (or equivalently an affine structure) on $G$ which is left invariant.  It is given on left invariant vector fields by the formula
             \begin{equation}\label{eq:LSA-symplectic} \omega^+ (\nabla_{x^+}y^+, z^+):= -\omega^+ (y^+, [x^+,z^+])                \end{equation}                which naturally extends  to a left invariant connection $\nabla$ on $G$.
The connection (\ref{eq:LSA-symplectic}) then induces an LSA structure on $\G,$ via the formula (\ref{eq:LSA-of-connection}).

\begin{example}\label{eg:aff(R)} The simplest example of a Frobenius Lie algebra is the 2-dimensional Lie algebra $\EuFrak{aff}(\mathbb K)$ of the group of affine motions of $\mathbb K.$ It has basis $(e_1,e_2)$  with Lie bracket $[e_1,e_2]=e_2.$  Let  $(e_1^*,e_2^*)$ stand for the corresponding dual basis. Set $\omega(x,y)=-e^*_2([x,y]).$ Then $\omega=\partial e^*_2 =-e_1^*\wedge e_2^*.$
The LSA is given by $e_1^2=-e_1,$ $e_1e_2=0,$ $e_2e_1=-e_2,$ $e_2^2=0.$  We also have, $q(e_1)=-e^*_2,$  $q(e_2)=e^*_1.$ The linear map defined by $J(e_1)=e_2,$ $J(e_2)=-e_1,$ together with $\omega,$ give rise to a left invariant K\"ahler structure on any Lie group with Lie algebra $\EuFrak{aff}(\mathbb K).$
\end{example}
 More generally,  a triplet $(\G,\omega^+_\epsilon,J^+_\epsilon)$  is a called a j-algebra,  if $(\G,\omega^+_\epsilon)$  is a Frobenius Lie algebra, with  $\omega^+_\epsilon(x,y) = -\alpha([x,y]),$ for all $x,y\in\G,$  and $(G,\omega^+,J^+)$ is a K\"ahler group, for some Lie group $G$ with Lie algebra $\G$ (see e.g. \cite{Li-Me})    j-algebras  are a key tool in  the study of  homogeneous K\"ahler Manifolds and in particular homogeneous bounded domains (see e.g. \cite{Do-Na},  \cite{Gindikin-Sapiro-Vinberg},  \cite{Shapiro69}).

   Let  $p\ge 1,$ $ n\ge 2$ be integers such that $p$ divides $n$. 
Consider the Lie algebra $\EuFrak{gl}(n,\mathbb K)$  of $n\times n$ real matrices,
 where the Lie bracket of two matrices $M,$ $N,$ is given, as usual, by  their commutator $[M,N]=MN-NM.$ 
 Let $\EuFrak{gl}(n,\mathbb K)$ act on the space $M_{n,p}$ of $n\times p$ real matrices by ordinary 
left multiplication of matrices. The resulting semi-direct product  
 $\EuFrak{gl}(n,\mathbb K)\ltimes M_{n,p}$  is a (non-solvable) Frobenius Lie algebra 
(see \cite{Rais}), with the Lie algebra of the group  $A\EuFrak{ff}(n,\mathbb K)$ of affine motions
 of $\mathbb K^n$  as a particular case, when $p=1$ (see also \cite{bordemannmedinaouadfel}). 
Frobenius Lie algebras of dimension less than or equal to $6,$ are extensively studied by several authors (see a list of those in dimension $4$ in \cite{Diatta2000}).

\section{Some properties of principal, right and nil elements}\label{chap:properties_of_principal_and_rightnil}

\begin{lemma}\label{properties1} Let $\G$ be a Lie algebra endowed with  an LSA.
 (1) Every  right-nil $\tilde x_0$ is an invariant or, equivalently, a $0$-cocycle (a  $0$-covariant element) for the representation $\mathcal L.$ That is, $0=\partial \tilde x_0 (y)=\mathcal L_y(\tilde x_0), $ for all $y\in\G.$
 \\
 (2) Let $y_0$ be a right-unit, then the following hold true:
(a) $-1$ and $0$ are always eigenvalues of the adjoint operator $ad_{y_0}$ of $y_0;$
(b) every right-nil element is an eigenvector of $ad_{y_0}$   corresponding to the eigenvalue
 $-1.$ In other words, right-nil elements are common eigenvectors of adjoint operators of all right-units,
 with $-1$ as the corresponding common eigenvalue;
(c) a vector $x$ is an eigenvector for $ad_{y_0}$ with corresponding eigenvalue $\lambda,$ if and only if $x$ is also an eigenvector for $\mathcal L_{y_0}$ with corresponding eigenvalue $1+\lambda.$
\end{lemma}
\proof (1) This  readily follows from the definition. If $\tilde x_0$ is a right-nil vector, then for every $y\in\G$ we have $0=y\tilde x_0=\mathcal L_{y}(\tilde x_0)=:\partial \tilde x_0 (y).$
(2)-(a) and (b): First, $0$ is obviously an eigenvalue of any adjoint operator.
Now if $y_0$ is a right-unit, the relation $[y_0,x]=y_0x-x$ implies in particular, that, if $x$ is a right-nil vector, then $ad_{y_0}(x) = -x.$
 (2)-(c) This is seen by rewriting the relation $[y_0,x] =y_0 x -x $ as   $(1+ad_{y_0})(x)=\mathcal L_{y_0} (x).$

 \qed

From Lemma \ref{properties1},  the kernel $\ker({\mathcal L}_{y_0})$ of ${\mathcal L}_{y_0}$
  is also the eigenspace $\mathcal G_{-1}$ of $ad_{y_0}$ corresponding to the eigenvalue $-1.$
 Hence $\ker({\mathcal L}_{y_0})$ is never  empty as it contains a subspace of dimension
$\dim(\mathcal G)-\dim([\mathcal G,\mathcal G])\neq 0.$ Recall that, from \cite{Hel79},
 if a Lie algebra possesses an LSA, then  $\dim(\mathcal G)\neq\dim[\mathcal G,\mathcal G].$  Lemma \ref{properties1},
(2) also implies that $ad_{y_0}$ cannot be nilpotent, if $y_0$ is a right-unit.

Now suppose $(\G,\omega=\partial\alpha)$ is a Frobenius Lie algebra. Consider the LSA structure as in (\ref{eq:LSA-symplectic}) and the isomorphism  $q:\G\to\G^*$ of vector spaces given by (\ref{eq:q}).
There exists a unique vector $x_\alpha$ in $\mathcal G$ such that $q(x_\alpha)= \alpha$.

\begin{definition}\label{lem:principalelement}
The unique vector $x_\alpha$ in $\mathcal G$ such that $q(x_\alpha)= \alpha$ is referred to as the principal element corresponding to $\alpha.$ 
\end{definition}
Throughout this work, we will simply write $x_0$ instead of $x_\alpha.$

The principal element $x_0$ is a right-unit with respect to the LSA (\ref{eq:LSA-symplectic}), as
 can be seen by simply using the equality $\omega(\mathcal R_{x_0}(x),y)= -\omega(x_0,[x,y])=
 - \langle q(x_0), [x,y]\rangle=- \alpha( [x,y])=\omega(x,y),$  which is valid for all $x,y\in \G.$

Let $(G,\omega^+:=d\alpha^+)$ be a Frobenius Lie group with Lie algebra  $(\G,\partial\alpha).$ 
Then, the Lie derivative  $L_{x^+_o}$ along $x_o^+$ satisfies $L_{x^+_o}\omega^+=\omega^+,$
 that is, $x_o^+$ is  a Liouville vector field on $G.$  We also have $L_{x^+_o}\alpha^+=\alpha^+$
 and $L_{x^+_o}r^+=-r^+$, where $r^+$ is the left invariant Poisson tensor corresponding 
to the solution $r\in\wedge^2\G$ of the classical Yang-Baxter Equation, 
defined using the duality pairing and the map  $\tilde r= q^{-1}$, 
by $r(\beta, \gamma):=\langle\tilde r(\beta),\gamma\rangle,$ for all $\beta, \gamma\in \G^*$ 
(see e.g. \cite{di-me-cybe} or \cite{Dr83} for more on the classical Yang-Baxter Equation).

\blem\label{lem:isomorphy-rightnil-cohom1} Let $\G$ be a Frobenius Lie algebra and $q$ defined as above, $\beta_0 \in \G^*$ and $y_0 :=q^{-1} (\beta_0).$ Then, we have
\beq  \partial \beta_0(x,y):=- \beta_0([x,y])=\omega(xy_0,y), \eeq
which implies that $\mathcal R_{y_0}(x)=xy_0=0$ for all $x\in \G$ if and only if the 1-form $\beta_0$ is closed. In other words, $q$ ismorphically maps right-nil elements to linear closed 1-forms on $\G.$
\elem
\proof Obviously, we have
$ \omega(xy_0,y)= - \omega(y_0,[x,y])= - \beta_0([x,y]) =\partial \beta_0(x,y).$
\qed

Let us remark that the space of closed linear 1-forms on $\G$ is isomorphic to the first cohomology space  $H^1(\G,\mathbb K)$  for the trivial action of $\G$ on $\mathbb K.$ We have

\begin{corollary} \label{lem:isomorphy-rightnil-cohom2} Let $(\G,\omega=\partial\alpha)$ be a Frobenius Lie algebra with principal element $x_0.$ Then every right-unit $y_0$ for the LSA given by (\ref{eq:LSA-symplectic}) is of the form $y_0=x_0+\tilde x_0,$ where $\tilde x_0$ is some right-nil element for the same LSA. In other words, the set of right-units is the affine space $x_0+q^{-1}(H^1(\G,\mathbb K)).$
\end{corollary}
\proof A vector $y_0$ is a right-unit on $\G,$ if and only if  for every  $x\in\G,$ we have $x(y_0-x_0)=0,$ that is, if and only if $\tilde y_0:=y_0-x_0$ is a right-nil vector in $\G.$ Thus, from Lemma \ref{lem:isomorphy-rightnil-cohom1} the set of right-units is $x_0+q^{-1}(H^1(\G,\mathbb K)).$
\qed

\begin{lemma}
Let $(G,\omega^+)$ be a quasi-Frobenius Lie group and   $(\mathcal G,\omega)$ the corresponding Quasi-Frobenius Lie algebra. The following are equivalent:
(a) the associated connection defined in (\ref{eq:LSA-symplectic}) has a radiant left invariant vector field or, equivalently, the induced LSA on $\mathcal G$ has a right-unit; \\
(b) $(G,\omega^+)$ is a Frobenius Lie group, i.e.  there exists $\alpha\in\mathcal G^*$ such that $\omega(x,y)=-\alpha([x,y]),$ for all $x,y\in\mathcal G.$\\
\end{lemma}
\proof
If $(G,\omega^+=d\alpha^+)$  is Frobenius, then $x_0:=q^{-1}(\alpha)$ is a right-unit for the LSA. That is, $x_0^+$ is a radiant left invariant vector field, for the connection $\nabla,$ defined by (\ref{eq:LSA-symplectic}).
Conversely, given a quasi-Frobenius Lie algebra, the cocyle $q:\mathcal G\to \mathcal G^*$ is equivariant with respect to the actions
 $\mathcal L$ of $\mathcal G$ on $\mathcal G$ via left multiplication and the coadjoint action on $\mathcal G^*.$ That is:
$q(xy)=ad^*_xq(y)$ for all $x,y\in\mathcal G.$  Now, suppose there is a right-unit $x_0,$ then
for every $x\in \mathcal G,$ we have
$q(x)=q(xx_0)=ad^*_xq(x_0)=ad_x^*\alpha,$ where $\alpha=q(x_0)\in\mathcal G^*$ satisfies
$\partial \alpha(x,y) = -\alpha([x,y])= - \langle q(x_0),[x,y]\rangle =- \langle q(x_0)\circ ad_x,y\rangle
=\langle ad_x^*q(x_0),y\rangle=\langle q(x),y \rangle=\omega(x,y).$ \qed

\begin{definition}
Let $\phi: \G\to\G$ be a linear operator (linear map). We say that $\phi$ is infinitesimally conformal with respect to $\omega$ if there exists a scalar number $\lambda(\phi)\in \mathbb K$ such that, for every $x,y\in\G,$ we have $\omega(\phi(x),y)+\omega (x, \phi(y))=\lambda(\phi)\omega(x,y).$
 If  more over, $\lambda(\phi)=0,$ we say that $\phi$ is infinitesimally symplectic with respect to $\omega.$ In the case where $\phi=ad_{v_0},$ for some vector $v_0\in\G,$ then we call $v_0$ a conformal vector.
\end{definition}
Symplectic vectors are those $v_0$ for which $\lambda(ad_{v_0})=0.$
 The linear map $\lambda : \EuFrak {gl}(\mathcal G) \to \mathbb K$, $\phi \mapsto \lambda(\phi)$
 is a Lie algebra homomorphism, that is, it is a closed linear form on $\G.$ Thus, the space $\ker (\lambda)=:\mathcal S$ of infinitesimal symplectic linear maps, is a Lie ideal of the Lie algebra $\mathcal C$ of infinitesimally conformal maps. More precisely, we have the following.
\blem \label{lem:conformal-symlectic-maps} In a quasi-Frobenius Lie algebra, the Lie algebras $\mathcal S$ of infinitesimally symplectic linear maps and $\mathcal C$ of infinitesimally conformal linear maps satisfy $[\mathcal C,\mathcal C]\subset \mathcal S.$
 The Lie subalgebras $\mathcal S_v$ of symplectic vectors and $\mathcal C_v$ of conformal vectors satisfy $[\mathcal C_v,\mathcal C_v]\subset \mathcal S_v.$
\elem
\proof As $\mathcal C$ and $\mathcal C_v$ are Lie algebras and trace($[\phi_1,\phi_2])=0,$ for any linear operators $\phi_1,~\phi_2,$ the proof readily derives from Proposition \ref{prop:trace_conformal-symlectic-maps} below. \qed

\bpro \label{prop:trace_conformal-symlectic-maps} Let $(\G,\omega)$ be a quasi-Frobenius Lie  algebra and $\phi: \G\to \G$  a linear map which is infinitesimally conformal with respect to $\omega.$ Then $trace(\phi)=\frac{1}{2}\lambda(\phi)\dim(\G).$ \epro
\proof
Let $(\G,\omega)$ be a quasi-Frobenius Lie  algebra of dimension $2n$, chose a Darboux basis $e_1,...,e_{2n}$ of $\G,$ so that $\omega(e_i,e_j)=\omega(e_{n+i},e_{n+j})=0$ and $\omega (e_i,e_{n+j})=\delta_{i,j},$ for all $i,j=1,2,\ldots,n$ where $\delta$ is the Kronecker delta symbol. Now if we write
 $\phi(e_l)=\displaystyle \sum_{k=1}^{2n}\phi_{k,l}e_k,$  with $\phi_{k,l}\in\mathbb K$  and $l=1,2,\ldots,2n.$ Then, the equality $\omega(\phi(e_i),e_{n+i})+\omega (e_i, \phi(e_{n+i}))=\lambda(\phi)\omega(e_i,e_{n+i}),$ becomes
$(\phi_{i,i}+\phi_{n+i,n+i})\omega (e_i, e_{n+i})=\lambda(\phi)\omega(e_i,e_{n+i}).$ Hence $trace(\phi)=\displaystyle \sum_{i=1}^n(\phi_{i,i}+\phi_{n+i,n+i})$$ = \sum_{i=1}^n \lambda(\phi)$$=n\lambda(\phi).$ \qed

On the manifold level (in the case at hand, quasi-Frobenius Lie groups) and for conformal maps, say $\Phi$, Proposition \ref{prop:trace_conformal-symlectic-maps} will simply translate into a similar relationship between the determinant of $\Phi$, $\lambda(\Phi)$ and the dimension of the manifold.

For a conformal vector $v_0$, we simply write $\lambda(v_0)$ instead of $\lambda(ad_{v_0}).$

\bpro \label{pro:conformalvectors} In a Frobenius Lie algebra $\G,$  right-unit and right-nil vectors are infinitesimally conformal vectors, with respect to $\omega.$ That is, there exists a Lie algebra homomorphism (i.e. a closed linear 1-form) $\lambda : \G\to \mathbb K$ such that $\omega([v_0,x],y)+\omega (x, [v_0,y])=\lambda(v_0)\omega(x,y),$ for every right-unit or right-nil vector $v_0$ and for every $x,y$ in $\G.$

More precisely right-nil elements are symplectic vectors, that is, they are in the kernel of $\lambda.$  Whereas right-unit vectors are solutions of the equation $\lambda(x)=-1,$ that is, every right-unit $v_0$ satisfies, for all $x,y\in\G,$
\beqn \omega([v_0,x],y)+\omega(x,[v_0,y])=-\omega(x,y).\label{rightunits-conformal}
\eeqn
\epro
\proof Suppose $v_0$ is a right-unit, then for every $x,y\in\G,$ we have $[v_0,x]=v_0x-x$ and $\omega([v_0,x],y)=\omega(v_0x,y)-\omega(x,y)=-\omega(x,[v_0,y])-\omega(x,y),$  or  equivalently $\omega([v_0,x],y)+\omega(x,[v_0,y])=-\omega(x,y).$
 Now if $v_0$ is a right-nil, then for every $x,y\in\G,$ the following holds  $[v_0,x]=v_0x.$ Thus  $\omega([v_0,x],y)=\omega(v_0x,y)=-\omega(x,[v_0,y]),$  hence $\omega([v_0,x],y)+\omega(x,[v_0,y])=0.$
\qed

In any Frobenius Lie group $(G,\omega^+)$ with Frobenius Lie algebra $(\G,\omega)$, the left invariant vector field $v_0^+$, corresponding to $v_0$ as in Proposition \ref{pro:conformalvectors}, satisfies $L_{v_0^+}\omega=\lambda(v_0)\omega^+.$
As a direct corollary of Proposition \ref{pro:conformalvectors}, the following Proposition \ref{prop:trace_principal_elements} generalizes a result in Theorem 3.3 of Ooms \cite{ooms80} on principal elements.

\bpro \label{prop:trace_principal_elements} In a Frobenius Lie algebra $\G$ of dimension $2n$, a right-unit vector  satisfies trace($ad_{v_0}$)$=-n$ and hence never lies in the derived ideal $[\G,\G].$ \epro
\proof
Suppose $v_0$ is a right-unit of $\G.$ From Proposition \ref{prop:trace_conformal-symlectic-maps}, trace($ad_{v_0}$)$=n\lambda(v_0)=-n.$ Hence $v_0$ is not in $[\G,\G],$ as trace$(ad_{[x,y]})=trace([ad_x,ad_y])=0,$ for all $x,y\in\G.$
\qed

The following is a corollary of Proposition \ref{prop:trace_principal_elements}.
\bpro A Frobenius Lie algebra cannot be unimodular. \epro
\proof A Lie algebra $\G$ is unimodular if and only if it has a adjoint-invariant volume form (in other words, every Lie group with Lie algebra $\G$ has a bi-invariant volume form. See \cite{Li-Me},  p. 226) or equivalently, the adjoint operator $ad_x$ associated to every element $x\in\G$, must be traceless.
\qed
\brmq In a Frobenius Lie group, the connection defined in (\ref{eq:LSA-symplectic}), is never geodesically complete (i.e. there exists at least one geodesic that cannot be extended to all values of its affine parameter). Indeed, from  \cite{Li-Me} this connection is complete if and only if the Lie group is  unimodular and hence solvable. This geodesic incompleteness  is also directly deduced from \cite{Fri-Gold-Hir}, using the existence of a radiant vector field. \ermq

\brmq A Frobenius Lie algebra cannot be embedded as a subalgebra of some compact Lie algebra. Indeed, it is well known that if a quasi-Frobenius Lie algebra $\G$ is a subalgebra of a compact Lie algebra, then $\G$ must be Abelian (see e.g. \cite{Li-Me}).  \ermq

\bthm \cite{Gerstenhaber2009}Suppose that $\G$ is a Frobenius Lie algebra which is a saturated subalgebra of a
simple Lie algebra $\mathcal S$ over $\mathbb K=\mathbb C$, that is, $\G$  is not an ideal of any larger subalgebra  of $\mathcal S.$ If $\alpha$ is a Frobenius functional on $\G,$ then its principal element $x_0$  is semisimple.
\ethm
In Section \ref{chap:semisimplicityprincipalelements}, we will see that the saturation property is not a necessary condition.

\bthm \cite{Gerstenhaber2009} Suppose a Frobenius Lie algebra $\G$ is a saturated Lie subalgebra of $\EuFrak{sl}(n,\mathbb C).$ If $x_0$ is a principal element, then the eigenvalues of $ad_{x_0}$ are integers which constitute a single unbroken string. That is, if $i$ and $j$ are eigenvalues with $i < j,$ then any integer $k$ with $i < k < j,$ is also an eigenvalue.
\ethm

\section{Frobenius Lie algebras are subalgebras of $\EuFrak{sl}(2n+1,\mathbb K)$}

The Following implies in particular that some of the results of M. Gerstenhaber and A. Giaquinto in
\cite{Gerstenhaber2009} on Frobenius subalgebras of
 $\EuFrak{sl}(m,\mathbb K)$ are actually valid for some more general
Frobenius Lie algebras. It also implies that the classification of the solutions
of the Classical Yang-Baxter Equation in $\EuFrak{sl}(2n+1,\mathbb K)$ alone (see e.g.
\cite{belavindrinfeld82, belavindrinfeld83, belavindrinfeld98}), is a route to classifying Frobenius and more generally quasi-Frobenius Lie algebras.

However, if $\mathcal G$ is a Frobenius Lie algebra, then $m=2n+1$ as in Theorem \ref{affineembedding}, is not the lowest integer for which   $\mathcal G$ is a subalgebra of $\EuFrak{sl}(m,\mathbb K).$
Indeed, for reasons related to the dimension of Cartan subalgebras, the embedded subalgebra $\phi(\mathcal G)$ does not contain any Cartan subalgebra of  $\EuFrak{sl}(2n+1,\mathbb K),$ unlike some of the saturated subalgebras in \cite{Gerstenhaber2009}. 

Now the question as to what obstructions prevent a Frobenius Lie algebra from being embedded as a
 subalgebra of a special kind, e.g. a saturated or seaweed, of some $\EuFrak{sl}(m,\mathbb K),$ is still open.

\begin{theorem}\label{affineembedding} Let $\G$ be a Lie algebra of dimension $p,$  endowed with an LSA structure, then:\\
(1) $\G$ is isomorphic to a Lie subalgebra of $\EuFrak{sl}(p+1,\mathbb K).$  In particular, quasi-Frobenius Lie algebras of dimension $2n$, can  all be seen as subalgebras of $\EuFrak{sl}(2n+1,\mathbb K),$
\\
(2)
suppose in addition, that in $\G$ every derivation is an inner derivation. Then every element of $\G$ admits a Jordan decomposition in $\G.$
\end{theorem}
\proof
 (1) From Sect. \ref{chap:properties_of_principal_and_rightnil}, a product 
$(x,y) \to xy =:\mathcal L_x(y)$ on $\G,$ is an LSA if and only if $\mathcal L: \mathcal G \to \EuFrak{gl} (\mathcal G)$ is a representation of the Lie algebra $\mathcal G$ and $id_{\G}:\mathcal G \to \mathcal G $ is 1-cocycle of $\mathcal L$, so that
\begin{eqnarray} \rho: \mathcal G \to \EuFrak{aff}(\mathcal G):=\EuFrak{gl} (\mathcal G)\ltimes \mathcal G, ~~  x~\mapsto ~\mathcal L_x+x
\end{eqnarray}
is a faithful representation of $\mathcal G$ by affine transformations of the vector space $\mathbb K^{p}$  underlying the Lie algebra $\mathcal G.$ Hence $\mathcal G$ is isomorphic to the subalgebra $Im(\rho)$ of $\EuFrak{aff}(\mathcal G).$
On the other hand, for any Lie subalgebra $\mathcal A$  of the Lie algebra of linear operators of a finite dimensional vector space,  the trace $Tr: \mathcal A \to \mathbb K,$~ $L\mapsto Tr(L)$ is a Lie algebra homomorphism.
Thus $\phi : \G\to \EuFrak{sl}(p+1,\mathbb K)$ defined by
\begin{eqnarray}\label{eq:phi}
\phi(x) = \begin{pmatrix} \mathcal L_x & x\\  0 & -Tr(\mathcal L_x) \end{pmatrix},\end{eqnarray}
is an injective Lie algebra homomorphism, hence $\phi(\G)$  is a subalgebra of $\EuFrak{sl}(p+1,\mathbb K)$ isomorphic to $\G.$
The proof of part (2) of Theorem \ref{affineembedding}, is a consequence of the following key lemma.        \qed
\begin{lemma}\label{lem:semisimple-jordan-cohom}
If $\G$  is a subalgebra of a semisimple Lie algebra  and  every derivation of  $\G$ is an inner derivation, then all elements admit a Jordan-Chevalley decomposition in $\G.$
\end{lemma}
\proof
The proof is similar to that of Theorem \ref{thm:semisimplecohomology} below.
  \qed
\begin{corollary}
For every Lie algebra $\G$ having trivial center $Z(\G)=0$ and only  inner derivations, all elements admit a Jordan-Chevalley decomposition in $\G.$
\end{corollary}
\proof It suffices to embed $\G$ in some $\EuFrak{sl}(m,\mathbb K)$  and then apply Lemma \ref{lem:semisimple-jordan-cohom} above.
\qed

\section{Semisimplicity of the Principal elements}\label{chap:semisimplicityprincipalelements}

The saturation condition is not a necessary condition for the semisimplicity of principal elements of Frobenius Lie algebras which are subalgebras of semisimple Lie algebras. Indeed, let us start  from a Frobenius Lie algebra $\G_1$ which
is a  saturated Lie subalgebra of a simple Lie algebra, say $\EuFrak{sl}(n,\mathbb K).$ Hence from \cite{Gerstenhaber2009} every principal element
corresponding to some Frobenius functional on $\G_1$, is semisimple. Now consider the direct sum  $\G_1\oplus \G_2=:\G$ (as Lie algebras)
of $\G_1$ and some other Frobenius Lie algebra $\G_2,$ so that $\G_1$  and $\G_2$ are both ideals of $\G.$
 As $\G$ is Frobenius, then from Theorem \ref{affineembedding}, $\G$ is a Lie subalgebra of some $\EuFrak{sl}(m,\mathbb K)$ for which, $\G_1$ is no longer saturated.

\begin{theorem} If in a Frobenius Lie algebra $\G$  every derivation  is an inner derivation,
then every embedding of $\G,$ in any other Lie algebra $\mathcal B,$ is either saturated or an ideal of a direct sum $\G\oplus\tilde \G$ of two subalgebras $\G$ and $\tilde \G,$ of $\mathcal B.$
 \end{theorem}
\proof The proof is straightforward. If $\G$ is an ideal of some Lie algebra that may be taken as a semi-direct $\G\rtimes \bar \G$ of $\G$ and another subalgebra $\bar \G,$ then
as  every derivation  is an inner derivation,   the action of $\bar \G$ on $\G$  induced by this semi-direct product, is done by inner derivations of $\G$ and hence
$\G\rtimes \bar \G$ is isomorphic to a direct sum of two subalgebras (ideals of this sum), one of which being $\G$. \qed

\begin{theorem}\label{thm:semisimplecohomology} If in a Frobenius Lie algebra $\G$  every derivation  is an inner derivation,
then the adjoint operator of every principal element is semisimple, at least for $\mathbb K=\mathbb C.$
 \end{theorem}
\proof We apply a reasoning akin to that in \cite{Gerstenhaber2009}. It uses the Jordan-Chevalley decomposition, which often appears to be a powerful tool, allowing to solve questions complicated in appearance, with very simple algebraic arguments. First, using Theorem \ref{affineembedding}, let us embed $\G$ in $\EuFrak{sl}(n+1,\mathbb K)$ and pick a principal functional $\alpha$ with principal element $x_0.$ As above, set $\omega:=\partial \alpha,$ that is, $\omega(x,y):=-\alpha([x,y]).$
Regard $x_0$ as an element of $\EuFrak{sl}(n+1,\mathbb K)$ and let  $x_0 = x_s + x_n$ be the Jordan-Chevalley decomposition  of $x_0$ into a semisimple and a nilpotent parts $x_s$ and $x_n,$ respectively. This decomposition holds true in every linear representation of $\G.$ In particular,  $ad_{x_0} = ad_{x_s} +ad_ {x_n}$ is the Jordan-Chevalley decomposition of $ad_{x_0},$ where $ad_{x_s}$ is semisimple and $ad_ {x_n}$ nilpotent, with $[ad_{x_s},ad_ {x_n}]=ad_{[x_s,x_n]}=0.$ Let $p_s$ and $p_n$ be two polynomials with one variable $T$, each without constant term, such that $p_s(ad_{x_0})=ad_{x_s}$ and  $p_n(ad_{x_0})=ad_{x_n}.$ Then, as polynomials in $ad_{x_0},$ the derivations $ad_{x_s}$ and $ad_{x_n}$ both preserve $\G,$ that is,  $ad_{x_s}(\G)\subset \G$ and $ad_{x_n}(\G)\subset \G.$ The hypothesis that every derivation of $\G$ is inner, insures that there exist $x_s'$ and $x_n'$ in $\G$ such that, the restrictions of $ad_{x_s}$ and $ad_{x_n}$ to $\G$ satisfy $ad_{x_s}=ad_{x_s'},$ $ad_{x_n}=ad_{x_n'}.$ For the rest of this proof, we can simply write $x_s$ and $x_n$  for  $x_s'$ and $x_n',$ without lost of generality. Using the identity $\alpha\circ ad_{x_0}= - \alpha,$  we get $\alpha\circ ad_{x_n}=p_{n}(-1)\alpha$ and $\alpha\circ ad_{x_s}=p_{s}(-1)\alpha.$
For every $y\in\G,$ we have $\langle\alpha\circ ad_{x_n},y\rangle= - \omega(x_n,y)$ and $ p_{n}(-1)\langle \alpha,y \rangle=p_{n}(-1)\omega(x_0,y).$
This is equivalent to the equality $x_n=-p_{n}(-1)x_0.$ But $ad_{x_0}$ has nonzero eigenvalues, such as $-1,$ moreover, from Proposition \ref{prop:trace_principal_elements}, trace($ad_{x_0}) =- \frac{1}{2}\dim (\G),$  which is not possible for the nilpotent operator $ad_{x_n}.$ So we are left with the only possibility $p_n(1)=0$ and  hence $ad_{x_n}=0.$ Or, equivalently, $x_n=0,$ given that the center of $\G$ is trivial. Hence $x_0$ is semisimple. \qed

Remark that every noncompact semisimple Lie algebra ($\mathbb K=\R$ or $\C$) contains a subalgebra which is a Frobenius Lie algebra all of whose principal elements are semisimple. Indeed, the root decomposition implies that any noncompact semisimple Lie algebra contains a copy of $\EuFrak{sl}(2,\mathbb K)$ which in turn, contains $\EuFrak{aff}(\mathbb K)$ as a subalgebra (see Theorem \ref{thm:affineliealgebra} ).

\section{The Lie algebra of the group of affine motions of the Euclidean space}

 For each  $n\ge 1,$ we define $\G_{n}$  as being the Lie algebra
$\G_{n}:=\EuFrak{aff}(n,\mathbb K)$  of the group of affine motions of $\mathbb K^n,$
with $\mathbb K=\R$ or $\C.$   The following result, stated in  Chapter 5 of \cite{Diatta2000},  
 which we reproduce here with a detailed proof, provides a family parametrized
by $n$ of examples of Lie algebras $\G_{n}$ , which satisfy the conditions of
Theorem \ref{thm:semisimplecohomology} above.  On its own, this is an interesting
result useful for many other applications, some of which will be explored in
 some upcoming work. When $n\ge 2,$ every $\G_n$ is nonsolvable,
whereas $\G_1$ is always solvable.
 Recall  that  affine motions are maps $\phi: \mathbb K^n\to \mathbb K^n$ of
 the form $\phi(v)=\psi(v)+v_0,$ for every $v\in\mathbb K^n,$ where $v_0$
is some vector in $\mathbb K^n$ and  $\psi: \mathbb K^n\to \mathbb K^n$ some invertible linear map.
 The set of such affine motions, endowed with the composition law of operators of
 the vector space $\mathbb K^n,$ is a Lie group whose Lie algebra $\EuFrak{aff}(n,\mathbb K)$
  is obtained by dropping the invertibility condition on the linear part $\psi.$
It is the semi-direct product  $\G_n:=\EuFrak{gl}(n,\mathbb K)\ltimes \mathbb K^n$
where the Lie algebra  $\EuFrak{gl}(n,\mathbb K),$ of  linear maps $\psi: \mathbb K^n\to \mathbb K^n,$
 acts (on the left)
 on the vector space  $\mathbb K^n$ in the ordinary way (linear maps applied to vectors).

 \begin{theorem}\label{thm:affineliealgebra} In the  Lie algebra
$\G_{n}:=\EuFrak{aff}(n,\mathbb K)$  of the group of affine motions of $\mathbb K^n,$ every derivation  is an inner derivation.
 \end{theorem}
\proof We use a direct calculation approach.  First, remark that the derived ideal of
 $\G_{n}$  is $[\G_{n},\G_{n}]=\EuFrak{sl}(n,\mathbb K)\ltimes \mathbb K^n$ and $\G_{n}$ is the semidirect product $\G_{n}=\mathbb K \mathbb I_{\mathbb K^n}\ltimes [\G_{n},\G_{n}]$ of  the line $\mathbb K \mathbb I_{\mathbb K^n}$ generated  by the identity map $x\in\mathbb K^n\mapsto x=:\mathbb I_{\mathbb K^n}(x)$ and  $[\G_{n},\G_{n}].$ We write $a_0:=\mathbb I_{\mathbb K^n}.$ Using the vector space decomposition $\G_{n}=\mathbb K a_0\oplus \EuFrak{sl}(n,\mathbb K)\oplus \mathbb K^n,$  an element of $\G_{n}$ is of the form $v=ka_0+a+x,$ where $k\in\mathbb K,$ $a\in \EuFrak{sl}(n,\mathbb K),$  and $x\in \mathbb K^n.$ The Lie bracket has the following properties: $[a_0,a]=0,$ $[a,b]\in \EuFrak{sl}(n,\mathbb K),$ $[a_0,x]=x,$ $[a,x]=ax\in \mathbb K^n$ and $[x,y]=0,$ for all $a,b\in \EuFrak{sl}(n,\mathbb K),$ $x,y\in \mathbb K^n.$
A derivation $D$ of $\G_{n}$ can be written using its components, as follows. For every
$v\in \G_{n},$ we have $D(v)=k_D(v) a_0 + D_1(v) + D_2(v),$ where $k_D: \G_{n}\to \mathbb K,$ is a linear closed 1-form on $\G_{n},$ and $D_1: \G_{n}\to \EuFrak{sl}(n,\mathbb K),$ ~
$D_2: \G_{n}\to \mathbb K^n$ are linear.  The reason why $k_D$ is closed, is the fact that the derived ideal is a characteristic ideal, i.e., it is preserved by every derivation.
If $a\in \EuFrak{sl}(n,\mathbb K),$ we have
$0=D([a_0,a])= [D(a_0),a]+[a_0,D(a)]$$=[D_1(a_0)+D_2(a_0),a]+ [a_0,D_2(a)].$ We end up having
\beqn\label{eq:der1aff}0=[D_1(a_0),a]-aD_2(a_0)+D_2(a)\eeqn
In Eq. (\ref{eq:der1aff}), the vector $[D_1(a_0),a]$ lies in $\EuFrak{sl}(n,\mathbb K)$ whereas $-aD_2(a_0)+D_2(a)$
is in $\mathbb K^n,$ for all $a \in \EuFrak{sl}(n,\mathbb K).$ We then deduce that, on the one hand $[D_1(a_0),a]=0,$ for all  $a\in \EuFrak{sl}(n,\mathbb K),$
 that is,
 \beqn\label{eq:der2aff}D_1(a_0)=0,\eeqn
  and on the other hand,
\beqn\label{eq:der3aff}D_2(a)=aD_2(a_0)=[a,D_2(a_0)].\eeqn

For $x\in \mathbb K^{n},$  we have $D(x)=D([a_0,x])=[D(a_0),x]+[a_0,D_2(x)]$ $=k_D(a_0) x+D_2(x),$  which now reads $D_1(x)+D_2(x)=k_D(a_0) x+D_2(x) $ or, equivalently,
\beqn\label{eq:der4aff}
D_1(x)&=&0 \text{ for all~} x\in \mathbb K^n \text{ and } k_D(a_0)=0.
\eeqn
Eq. (\ref{eq:der4aff}) means that $D(x)=D_2(x)\in \mathbb K^n,$ for every $x\in \mathbb K^n$ and $D(a_0)=D_2(a_0)\in \mathbb K^n.$
So from now on, let us set $v_0=D(a_0)\in \mathbb K^n$  and denote by $\psi$ the restriction of $D$ to $\mathbb K^n.$
 So  $\psi$ is simply a linear map $\psi: \mathbb K^{n} \to \mathbb K^{n},$
with $\psi(x)=D(x)=D_2(x).$ Under these new notations,
 Eq. (\ref{eq:der3aff}) becomes $D_2(a)=av_0=[a,v_0],$
 whereas Eq. (\ref{eq:der4aff}) implies that $D(a_0)=v_0=[-v_0+\psi,a_0]$ and $D(x)=\psi (x)=[-v_0+\psi,x],$ for every $a\in \EuFrak{sl}(n,\mathbb K)$ and $x$ in $\mathbb K^n.$
Now we also have $\psi(a x)=D([a,x])$ $=[D(a),x]+[a,\psi(x)]$  $=D_1(a)x+a\psi(x),$
 that is, $[\psi,a] x=D_1(a)x,$ for all $a\in \EuFrak{sl}(n,\mathbb K)$ and $x\in \mathbb K^{n}.$ This can be written as
 \beqn\label{eq:der5aff} [\psi,a]=D_1(a)\in \EuFrak{sl}(n,\mathbb K), \text{~ for all~}
 a\in \EuFrak{sl}(n,\mathbb K).
\eeqn
So combining Eqs. (\ref{eq:der3aff}) and (\ref{eq:der5aff}), we get $D(a)=D_1(a)+D_2(a)=[\psi-v_0,a],$ for every $a\in \EuFrak{sl}(n,\mathbb K).$
 Summing up, for each derivation $D$ of $\G_n,$
there exist a vector $v_0\in \mathbb K^n$ and a linear map $\psi:~ \mathbb K^n\to \mathbb K^n$
 such that $D(v)=[-v_0+\psi,v]=ad_{V_D}(v),$ for all $v\in\G_{n},$ where $V_D:=-v_0+\psi$ is a element of $\G_n.$  \qed

\section{Examples}\label{chap:examples}
In this section, we provide a wide family of Frobenius solvable Lie algebras. We denote by $\delta_{ij}$ the Kronecker symbol, i.e. $\delta_{ij}=0$ if $i\neq j$ and $\delta_{ii}=1.$

  Let $\mathcal B= (e_{-1},e_0,e_1,\ldots,\ldots,e_{2n})$  be a basis of $\mathbb K^{2n+2}=\mathbb K^{2}\times\mathbb K^{2n}$ with corresponding dual basis  $(e_s^*)_{s},$ $s=-1,0,1,2,\ldots,2n.$ Consider a nonzero real number $k$ and a linear map $\xi:\mathbb K^{2n}\to \mathbb K^{2n},$  the coefficients $M_{s,t},$ $s,t=1,\cdots, 2n,$ of whose matrix $M$ in the basis $(e_1,\ldots,\ldots,e_{2n}),$  are real numbers satisfying
 $k\delta_{ij}=M_{j,i}+M_{n+i,n+j},$ ~ $M_{i,n+j}=M_{j,n+i}$ and $M_{n+i,j}=M_{n+j,i},$ for all $i,j=1,\cdots, n.$
Then the bracket in $\mathbb K^{2n+2}$ defined, for all $i,j=1,\cdots, n$ and $t=1,\cdots,2n,$ by $[e_{i},e_{n+j}]=\delta_{ij}e_0,$
 $[e_{-1},e_{0}]=ke_0,$ $[e_{-1},e_{t}]=\xi (e_t)=\displaystyle \sum_{s=1}^{2n}M_{s,t}e_s,$ is a Lie bracket and $e^*_0$ is a Frobenius functional on $\mathcal G_{k,\xi}:=(\mathbb K^{2n+2},[.,.]).$ More precisely, we have
$\partial e^*_{0}= \displaystyle \sum_{i=1}^ne_{n+i}^*\wedge e_{i}^* +ke_{0}^*\wedge e_{-1}^*,$  the corresponding principal element is $x_0=-\frac{1}{k}e_{-1}.$ In the basis $\mathcal B$,  the adjoint of the principal element reads $ad_{x_0}=-\frac{1}{k}\begin{pmatrix} 0&0&0\\ 0& k & 0\\
0 & 0 & M\end{pmatrix}.$ Thus, $ad_{x_0}$ is diagonalizable, triangularizable, ... if and only if, $M$ is. The conditions on the coefficients of the  matrix $M$ of $\xi$ are not very restrictive. 
Indeed, the set of such couples $(\xi,k)$ with $k\neq 0,$ is an open subset of the Lie algebra $\mathcal C=\EuFrak{csp}(n,\mathbb R)$ of infinitesimally conformal maps of $(\mathbb R^{2n},\Omega),$ where $\Omega$ is the canonical symplectic form $\Omega:= \displaystyle \sum_{i=1}^ne_{n+i}^*\wedge e_{i}^* .$ Namely, for every couple $(\xi',k'),$ where $\xi': \mathbb R^{2n}\to \mathbb R^{2n}$ is linear and $k'\in \mathbb R$ such that for every $x,y\in \mathbb R^{2n},$ we have $\Omega (\xi' (x),y)+ \Omega (x,\xi'(y) ) = k' \Omega (x,y),$  we define a Frobenius Lie algebra $\G_{k',\xi'}$ with underlying vector space $\mathbb K^{2n+2}$ as above, if and only if $k'\neq 0.$

Recall that $\EuFrak{csp}(n,\mathbb R)$  contains the Lie algebra $\EuFrak{sp}(n,\mathbb R)$ of infinitesimal symplectomorphisms as an ideal, as discussed in Section \ref{chap:properties_of_principal_and_rightnil}, where we have written $\mathcal C$ and $\mathcal S$ for the case of a general Quasi-Frobenius  Lie algebra, instead of $\EuFrak{csp}(\G)$ and $\EuFrak{sp}(\G)$, respectively.

The Lie algebra of derivations of each $\G_{k,\xi}$ is far from been isomorphic to $\G_{k,\xi}.$ 
Indeed, the simplest examples of non-inner derivations are the following. 
Every element $(k',\xi') $  of $\EuFrak{csp}(n,\mathbb R)$  where $\xi'$  
commutes with $\xi,$ defines a derivation  $D_{(k',\xi')} $  of $\G_{k,\xi}$, by 
 $D_{(k',\xi')} (e_{-1})=0,$ ~  $D_{(k',\xi')} (e_0)=k'e_0,$ ~ $D_{(k',\xi')} (e_t)=\xi'(e_t),$
 for $t=1,\cdots,2n.$ If  $\xi' $ is not of the form $\xi' =k''\psi\circ \xi\circ \psi^{-1} ,$ where 
  $\psi$ is a conformal map and $k''$ is a scalar,  then $D_{(k',\xi')} $   is not an inner derivation of   $\mathcal G_{k,\xi}.$

Let us consider the following specific example.  Take a $\xi$ as above, the coefficients of the corresponding $M$ being all zero, except the following  $M_{i,i}=k-M_{n+i,n+i}=k_i$ and $M_{n+i,i}=-M_{i,n+i}=k_i',$ so that Lie bracket in $\G_{k,\xi}$ is 
$[e_i,e_{n+i}]=e_0,$ $[e_{-1},e_i]=k_i e_i + k_i'e_{n+i},$
$[e_{-1},e_{n+i}]=- k_i'e_{i}+(k-k_i) e_{n+i}$ for $i=1,\ldots,n$    and $[e_{-1},e_0]=ke_0.$ If we take $k_i'=0,$ for all $i=1,\cdots,n$ then $ad_{x_0}$ is diagonal  $ad_{x_0}=diag(0,k,k_1,\ldots, k_n, k-k_1,\ldots,k-k_n).$  We can chose, for example, $k=1$ and $k_i=\pi^i$ so that $ad_{x_0}=diag(0,1,\pi,\ldots, \pi^n, 1-\pi,\ldots,1-\pi^n)$  and hence eigenvalues are all irrational, except two. In fact,  in the general case, we can invite any number  to be an eigenvalue. For instance, if we chose $k=1$ and $k_i,k_i'$ satisfying $k_i(1-k_i)+k_i'^2=-1,$ then the golden number $\frac{1+\sqrt{5}}{2}$  and $\frac{1-\sqrt{5}}{2}$   are repeated eigenvalues, as the characteristic polynomial will have $X^2-X-1$ as a multiple factor.

Now let us specialize in the case $n=1.$  
\\
(a) Take $\xi (e_1)=0,$ ~~~ $\xi (e_2)=e_2$   and $k=1,$ so that the Lie bracket is given by  
 $[e_{-1},e_0]=e_0,$ $[e_{-1},e_2]=e_2,$ $[e_1,e_2]=e_0.$  Thus we have $x_0=-e_{-1},$ and 
$ad_{x_0}=- diag(0,1,0,1).$
\\
(b) Take $\xi (e_1)=-e_1,$ ~~~ $\xi (e_2)=-e_1-e_2$   and $k=-2,$ so that the Lie bracket is given by
$[e_1,e_2]=e_0,$ $[e_{-1},e_0]=-2 e_0,$ $[e_{-1},e_1]=-e_1,$ $[e_{-1},e_2]=-e_1-e_2,$ $x_0=\frac{1}{2}e_{-1}$ and

 $ad_{x_0} = -\begin{pmatrix}0 & 0&0&0\\ 0 & 1& 0&0\\ 0 & 0&\frac{1}{2}&\frac{1}{2}\\ 0 & 0&0&\frac{1}{2}
\end{pmatrix}$
 cannot be diagonalized. Indeed its eigenvalues are $-1,$ $-\frac{1}{2},$ $-\frac{1}{2},$ $0$ and the eigenspace corresponding to the repeated eigenvalue $-\frac{1}{2}$ is 1-dimensional and spanned by $(0,0,1, 0).$ The eigenspaces corresponding to $-1$ and $0$ are spanned by $(0,1,0,0)$ and  $(1,0,0,0)$ respectively.
\\
(c) Consider $\xi (e_1)=\tilde ke_1-e_2,$ ~~~ $\xi (e_2)=e_1+\tilde k e_2$   and $k=2\tilde k,$ for some $\tilde k$ in $\mathbb R,$ so that the Lie bracket is given by
$[e_1,e_2]=e_0,$ $[e_{-1},e_0]=2\tilde k e_0,$ $[e_{-1},e_1]=\tilde ke_1-e_2,$ $[e_{-1},e_2]=e_1+\tilde k e_2,$ $x_0=-\frac{1}{2\tilde k}e_{-1}.$ Now $ad_{x_0}$  is diagonalizable for $\mathbb K=\mathbb C,$  but not for $\mathbb K=\mathbb R.$ Indeed it has  four distinct eigenvalues, two of which are complex: $0,1,\frac{-\tilde k+i}{2\tilde k}, -\frac{\tilde k+i}{2\tilde k}.$

\end{document}